\documentclass[11pt]{amsart}
\usepackage{latexsym,amsmath,amssymb}

\title[Maximal potentials and the spherical maximal function]{Maximal potentials, maximal singular integrals, and the spherical maximal function}

\author{Piotr Haj\l{}asz, Zhuomin Liu}

\address{Department of Mathematics, University of Pittsburgh, 301
  Thackeray Hall, Pittsburgh, PA 15260, USA, {\tt hajlasz@pitt.edu}}

\address{Department of Mathematics, University of Pittsburgh, 301
  Thackeray Hall, Pittsburgh, PA 15260, USA, {\tt liuzhuomin@hotmail.com}}

\thanks{P.H. was supported by NSF grant DMS-1161425.}

plus 6pt minus 12pt
\def\eps{\varepsilon}

\def\vi{\varphi}

\def\M{{\mathcal M}}
\def\H{{\mathcal H}}

\def\S{{\mathcal S}}
\def\I{{\mathcal I}}

\def\rss{{\vert_{_{_{\!\!-\!\!-\!}}}}}
\newtheorem{theorem}{Theorem}
\newtheorem{lemma}[theorem]{Lemma}

\newtheorem{proposition}[theorem]{Proposition}

\theoremstyle{definition}

\newcommand{\barint}{
\rule[.036in]{.12in}{.009in}\kern-.16in \displaystyle\int }

\newcommand{\barcal}{\mbox{$ \rule[.036in]{.11in}{.007in}\kern-.128in\int $}}

\newcommand{\bbbr}{\mathbb R}

\def\bbbc{{\mathchoice {\setbox0=\hbox{$\displaystyle\rm C$}\hbox{\hbox
to0pt{\kern0.4\wd0\vrule height0.9\ht0\hss}\box0}}
{\setbox0=\hbox{$\textstyle\rm C$}\hbox{\hbox
to0pt{\kern0.4\wd0\vrule height0.9\ht0\hss}\box0}}
{\setbox0=\hbox{$\scriptstyle\rm C$}\hbox{\hbox
to0pt{\kern0.4\wd0\vrule height0.9\ht0\hss}\box0}}
{\setbox0=\hbox{$\scriptscriptstyle\rm C$}\hbox{\hbox
to0pt{\kern0.4\wd0\vrule height0.9\ht0\hss}\box0}}}}

\def\bbbq{{\mathchoice {\setbox0=\hbox{$\displaystyle\rm Q$}\hbox{\raise
0.15\ht0\hbox to0pt{\kern0.4\wd0\vrule height0.8\ht0\hss}\box0}}
{\setbox0=\hbox{$\textstyle\rm Q$}\hbox{\raise 0.15\ht0\hbox
to0pt{\kern0.4\wd0\vrule height0.8\ht0\hss}\box0}}
{\setbox0=\hbox{$\scriptstyle\rm Q$}\hbox{\raise 0.15\ht0\hbox
to0pt{\kern0.4\wd0\vrule height0.7\ht0\hss}\box0}}
{\setbox0=\hbox{$\scriptscriptstyle\rm Q$}\hbox{\raise
0.15\ht0\hbox to0pt{\kern0.4\wd0\vrule height0.7\ht0\hss}\box0}}}}

\def\bbbz{{\mathchoice {\hbox{$\sf\textstyle Z\kern-0.4em Z$}}
{\hbox{$\sf\textstyle Z\kern-0.4em Z$}} {\hbox{$\sf\scriptstyle
Z\kern-0.3em Z$}} {\hbox{$\sf\scriptscriptstyle Z\kern-0.2em
Z$}}}}

\def\mvint_#1{\mathchoice
          {\mathop{\vrule width 6pt height 3 pt depth -2.5pt
                  \kern -8pt \intop}\nolimits_{\kern -3pt #1}}%
          {\mathop{\vrule width 5pt height 3 pt depth -2.6pt
                  \kern -6pt \intop}\nolimits_{#1}}%
          {\mathop{\vrule width 5pt height 3 pt depth -2.6pt
                  \kern -6pt \intop}\nolimits_{#1}}%
          {\mathop{\vrule width 5pt height 3 pt depth -2.6pt
                  \kern -6pt \intop}\nolimits_{#1}}}

\numberwithin{theorem}{section} \numberwithin{equation}{section}

\begin{document}



\sloppy

\begin{abstract}
We introduce a notion of maximal potentials and we prove that they form bounded
operators from $L^p$ to the homogeneous Sobolev space $\dot{W}^{1,p}$ for all 
$n/(n-1)<p<n$. We apply this result to the problem of
boundedness of the spherical maximal operator in Sobolev spaces.
\end{abstract}

\maketitle

\section{Introduction}

Let $\Omega$ be a bounded function defined on the unit sphere $S^{n-1}:=S^{n-1}(0,1)$ with 
zero mean value
\begin{equation}
\label{zero}
\int_{S^{n-1}}\Omega(x)\, d\sigma(x)=0,
\end{equation}
and let
\begin{equation}
\label{singular}
K(x)=\frac{\Omega(x/|x|)}{|x|^n},
\quad
x\neq 0
\end{equation}
be the associated homogeneous kernel of degree $-n$. 
It is well known and easy to prove that under these assumptions the
principal value ${\rm p.v.}\, K$ is a tempered distribution.
Moreover the convolution with ${\rm p.v.}\, K$ defines a singular integral
that is bounded in $L^p$ for all $1<p<\infty$.
Namely the {\em singular integral}
$$
T_\Omega f(x) =\left({\rm p.v.}\, K\right)*f(x)=\lim_{t\to 0}\int_{|x-z|\geq t}f(z)K(x-z)\, dz
$$
and also the {\em maximal singular integral}
$$
T_\Omega^*f(x)=\sup_{t>0}\left|\int_{|x-z|\geq t} f(z)K(x-z)\, dz\right|
$$
are bounded in $L^p(\bbbr^n)$, $1<p<n$. This follows from 
\cite[Theorems~4.2.7 and~4.2.11]{grafakos} and the fact that the kernel $K$ 
can be represented as a sum of an odd and an even one. Actually,
the result is true under much weaker assumptions about $\Omega$, see \cite{grafakos}.

In this note we investigate related results for convolutions with 
homogeneous kernels of
degree $-(n-1)$. Such operators have smoothing properties in the sense that they 
increase regularity of functions.
More precisely if 
$\Omega\in C^1(S^{n-1})$ and 
$$
\tilde{K}(x)=\frac{\Omega(x/|x|)}{|x|^{n-1}},
\quad
x\neq 0
$$
(we no longer require \eqref{zero}), then the operator
\begin{equation}
\label{potential}
A_\Omega f(x) = \int_{\bbbr^n}f(z)\tilde{K}(x-z)\, dz = \lim_{t\to 0} \int_{|x-z|\geq t} f(z)\tilde{K}(x-z)\, dz
\end{equation}
is bounded from $L^p$ to the {\em homogeneous Sobolev space} $\dot{W}^{1,p}$ for all $1<p<n$.
This is well known (see Proposition~\ref{potential2}).
Actually we do not have to pass to the limit as the integral converges for a.e. $x\in\bbbr^n$
by the Fractional Integration Theorem, \cite[Theorem~2.8.4]{ziemer}.

It is natural to consider the associated maximal operator
\begin{equation}
\label{maximal-potential}
A_\Omega^*f(x) = \sup_{t>0}\left|\int_{|x-z|\geq t} f(z)\tilde{K}(x-z)\, dz\right|.
\end{equation}
We will call it a {\em maximal potential}.
Such operators with $\tilde{K}(x)=x/|x|^{s+1}$ are called {\em maximal $s$-Riesz transforms}
and in our case $s=n-1$. These operators have been studied for example in \cite{ENV}, \cite{tolsa},
however, they were investigated from a different perspective than the one discussed in this paper.
In Section~\ref{cztery} we will show that the spherical maximal operator of a Sobolev function
can be represented as a maximal $(n-1)$-Riesz transform.

The question that we investigate in the paper is for what values of $p$, the operator 
$A_\Omega^*:L^p\to \dot{W}^{1,p}$ is bounded.

Clearly the operator $A_\Omega^*f$ is bounded by the {\em Riesz potential}
$$
A_\Omega^*f(x)\leq C \I_1|f|(x) := C\int_{\bbbr^n} \frac{|f(z)|}{|x-z|^{n-1}}\, dz,
$$
and the Riesz potential is bounded
$$
\I_1:L^p\to\dot{W}^{1,p} 
$$
for all for $1<p<n$.
Thus $A_\Omega^*f$ has the growth properties of
a function from $\dot{W}^{1,p}$. On the other hand the radius $t$ at which the supremum in \eqref{maximal-potential}
is attained may depend on $x$ in a very irregular way and this perhaps can result in the lack of smoothness of $A_\Omega^*f$.

In our main result (Theorem~\ref{main}) we will prove that the operator $A_\Omega^*:L^p\to \dot{W}^{1,p}$ is bounded
for all $n/(n-1)<p<n$. We do not know what happens when $1<p\leq n/(n-1)$, but
we will show that the positive answer to the questions about boundedness for $1<p\leq n/(n-1)$
would imply a positive answer to a question about boundedness of the spherical maximal function in the Sobolev space,
\cite{hajlaszl}, \cite{hajlaszo},  see Proposition~\ref{spherical}.

Now we can state our main result.
\begin{theorem}
\label{main}
Let $\Omega\in C^1(S^{n-1},\bbbr^m)$ and let
$$
\tilde{K}(x)=(\tilde{K}_1(x),\ldots,\tilde{K}_m(x))=\frac{\Omega(x/|x|)}{|x|^{n-1}},
\quad
x\neq 0.
$$
For $f=(f_1,\ldots,f_m)\in L^p(\bbbr^n,\bbbr^m)$ we define
\begin{eqnarray*}
A_\Omega^*f(x)
& = &
\sup_{t>0}\left|\int_{|x-z|\geq t} f(z)\cdot \tilde{K}(x-z)\, dz\right| \\
& = &
\sup_{t>0}\left|\sum_{i=1}^m\int_{|x-z|\geq t} f_i(z)\tilde{K}_i(x-z)\, dz\right|\, .
\end{eqnarray*}
Then the operator $A_\Omega^*:L^p\to \dot{W}^{1,p}$ is bounded for $n/(n-1)<p<n$. 
\end{theorem}
We could formulate the theorem under weaker assumptions than 
$\Omega\in C^1(S^{n-1},\bbbr^m)$, but we wanted to keep the
presentation simple. After all the $C^1$ case is the most interesting one.

Kinnunen, \cite{kinnunen}, proved that the Hardy-Littlewood maximal operator is bounded
in the Sobolev space $\M:W^{1,p}\to W^{1,p}$, $1<p<\infty$. This result initiated the search for the 
boundedness results for various maximal operators in Sobolev spaces.

It was proved by Bourgain \cite{bourgain} and Stein  \cite{stein}, 
(see also \cite{grafakos}, \cite{steinbook}) that the spherical maximal operator
$$
\S f(x) = \barint_{S^{n-1}(x,t)} |f(z)|\, d\sigma(z)
$$
is bounded in $L^p(\bbbr^n)$ when $n\geq 2$ and $p>n/(n-1)$.
Using this result it is easy to prove, see \cite{hajlaszl}, \cite{hajlaszo},  that actually the spherical maximal
operator is bounded in the Sobolev spaces $\S:W^{1,p}\to W^{1,p}$ for all $n/(n-1)<p<\infty$.
In particular $\S:W^{1,p}\to \dot{W}^{1,p}$ is bounded for $n/(n-1)<p<n$.
It was posted as an open problem in \cite{hajlaszl}, \cite{hajlaszo},  whether 
$\S:W^{1,p}\to \dot{W}^{1,p}$ is bounded for $1<p\leq n/(n-1)$. 
The case of the spherical maximal function has also been discussed in \cite{kinnunens}.
As a corollary of Theorem~\ref{main} we will prove
\begin{proposition}
\label{spherical}
If the claim of Theorem~\ref{main} is true for $1<p\leq n/(n-1)$, then the spherical maximal operator
$\S:W^{1,p}\to \dot{W}^{1,p}$ is bounded for all $1<p\leq n/(n-1)$.
\end{proposition}

The paper is organized as follows. In Section~\ref{sobolev} we recall basic facts about Sobolev spaces and potentials
and provide a proof that the operator \eqref{potential}
is bounded from $L^p$ to the homogeneous Sobolev space $\dot{W}^{1,p}$ for all $1<p<n$.
In Section~\ref{trzy} we prove Theorem~\ref{main} and in Section~\ref{cztery} we prove
Proposition~\ref{spherical}.

Notation is pretty standard. By $C$ we denote a generic positive constant that may have 
different values in different inequalities. The volume of the unit ball is denoted by $\omega_n$,
so the volume of the unit sphere is $n\omega_n$. Finally
$$
\barint_{S^{n-1}(x,t)}f(x)\, d\sigma(x) = \frac{1}{n\omega_n t^{n-1}} \int_{S^{n-1}(x,t)} f(x)\, d\sigma(x)
$$
stands for the integral average.


\section{Sobolev spaces and singular integrals}
\label{sobolev}

The {\em Sobolev space} $W^{1,p}(\bbbr^n)$, $1\leq p\leq\infty$ is defined as a class of all $f\in L^p$
whose distributional partial derivatives of first order also belong to $L^p$,
$\nabla f\in L^p$. $W^{1,p}$ is a Banach space with respect to the norm 
$\Vert f\Vert_{1,p}=\Vert f\Vert_p+\Vert \nabla f\Vert_p$. The homogeneous Sobolev space $\dot{W}^{1,p}$, $1<p<n$
consists of all functions $f\in L^{p^*}(\bbbr^n)$ such that $\nabla f\in L^p$, where $p^*=np/(n-p)$. 
It is a Banach space with respect to the norm $\Vert f\Vert_{\dot{W}^{1,p}}=\Vert f\Vert_{p^*}+\Vert\nabla f\Vert_p$.
It follows
from the Sobolev embedding theorem that $W^{1,p}\subset \dot{W}^{1,p}$ for $1<p<n$. 

If a function $f\in L^{p^*}$ is absolutely continuous on almost all lines parallel to coordinate axes, then $f$
has partial derivatives defined a.e. If in addition these partial derivatives belong to $L^p$, then
$f\in \dot{W}^{1,p}$, see \cite[Section~4.9]{EG}. We will use this characterization of $\dot{W}^{1,p}$ in the proof
of Theorem~\ref{main}.

Now we will show that certain potential type operators are bounded from $L^p$ to
$\dot{W}^{1,p}$, $1<p<n$. The following result seems to be well known, but we could not find a
direct reference, so we provide a proof.
In particular it implies that the Riesz potential $\I_1:L^p\to\dot{W}^{1,p}$, $1<p<n$ is bounded.

\begin{proposition}
\label{potential2}
If $\Omega\in C^1(S^{n-1})$ and
$$
\tilde{K}(x)=\frac{\Omega(x/|x|)}{|x|^{n-1}},
\quad
x\neq 0,
$$
then the operator
$$
A_\Omega f=\tilde{K}*f
$$
is bounded from $L^p$ to $\dot{W}^{1,p}$, $1<p<n$.
\end{proposition}
{\em Proof.}
First observe that
$$
|A_\Omega f(x)|\leq C\int_{\bbbr^n} \frac{|f(z)|}{|x-z|^{n-1}}\, dz.
$$
Since by the Fractional Integration Theorem, \cite[Theorem~2.8.4]{ziemer}, the {\em Riesz potential}
$$
\I_1g(x)=\int_{\bbbr^n}\frac{g(z)}{|x-z|^{n-1}}\, dz
$$
is bounded from $L^p$ to $L^{p^*}$, $1<p<n$, we conclude that
$A_\Omega:L^p\to L^{p^*}$ is bounded.

\begin{lemma}
\label{BH7}
Under the above assumptions
the pointwise gradient
$\nabla \tilde{K}(x)$ defined for $x\neq 0$ is homogeneous of degree $-n$ and
\begin{equation}
\label{BHeq7}
\int_{S^{n-1}} \nabla \tilde{K}(x)\, d\sigma(x)=0\, .
\end{equation}
Hence the condition \eqref{zero} is satisfied and thus
$$
{\rm p.v.}\, \nabla \tilde{K}\in \S_n'
$$
is a well defined tempered distribution, i.e. for each $1\leq j\leq n$
$$
{\rm p.v.}\, \frac{\partial \tilde{K}}{\partial x_j}\in\S_n'\, .
$$
Finally the distributional gradient $\nabla \tilde{K}$ satisfies
\begin{equation}
\label{BHeq8}
\nabla \tilde{K} =
c\delta_0+\, {\rm p.v.}\, \nabla \tilde{K}\, ,
\end{equation}
where $\delta_0$ is the Dirac distribution and
$$
c=\int_{S^{n-1}}\tilde{K}(x)\frac{x}{|x|}\, d\sigma(x)\, .
$$
In other words for $\vi\in\S_n$ and $1\leq j\leq n$
we have
$$
\frac{\partial \tilde{K}}{\partial x_j}[\vi]:=
-\int_{\bbbr^n} \tilde{K}(x)\, \frac{\partial\vi}{\partial x_j}(x)\, dx =
c_j\vi(0)+
\lim_{\eps\to 0}
\int_{|x|\geq\eps}
\frac{\partial \tilde{K}}{\partial x_j}(x)\vi(x)\, dx\, ,
$$
where
$$
c_j = \int_{S^{n-1}}\tilde{K}(x)\, \frac{x_j}{|x|}\, d\sigma(x)\, .
$$
\end{lemma}
{\em Proof.}
The fact that $\nabla \tilde{K}(x)$ is homogeneous of degree $-n$ is easy
and left to the reader.
For $0<t<r$ let $A(t,r)=\{x:\, t\leq |x|\leq r\}$.
The integration by parts yields
\begin{eqnarray*}
\int_{1\leq |x|\leq r}
\nabla \tilde{K}(x)\, dx
& = &
\int_{\partial A(1,r)} \tilde{K}(x)\, \vec{\nu}(x)\, d\sigma(x) \\
& = &
-\int_{|x|=1} \tilde{K}(x)\, \frac{x}{|x|}\, d\sigma(x) +
\int_{|x|=r} \tilde{K}(x)\, \frac{x}{|x|}\, d\sigma(x) =0\, .
\end{eqnarray*}
Indeed, the last two integrals are equal by a simple change of variables
and homogeneity of $\tilde{K}$. Thus the integral on the left hand side
equals $0$ independently of $r$. Hence its derivative with respect to $r$
is also equal zero.
$$
0=
\frac{d}{dr}\Big|_{r=1^+}
\int_{1\leq |x|\leq r} \nabla \tilde{K}(x)\, dx =
\int_{|x|=1} \nabla \tilde{K}(x)\, d\sigma(x)\, .
$$
This proves (\ref{BHeq7}). Therefore ${\rm p.v.}\, \nabla \tilde{K}\in\S_n'$
is a well defined tempered distribution. We are left with the proof that
the distributional gradient $\nabla \tilde{K}$ satisfies (\ref{BHeq8}).
Let $\vi\in\S_n$. We have
\begin{eqnarray*}
\lefteqn{\nabla \tilde{K}[\vi]
 :=
-\int_{\bbbr^n} \tilde{K}(x)\nabla\vi(x)\, dx}\\
& = &
\lim_{\eps\to 0}\lim_{R\to\infty}
-\int_{\eps\leq |x|\leq R} \tilde{K}(x)\nabla \vi(x)\, dx\\
& = &
\lim_{\eps\to 0}\lim_{R\to\infty}
\left(\int_{\eps\leq |x|\leq R}
\nabla \tilde{K}(x)\, \vi(x)\, dx -
\int_{\partial A(\eps,R)}
\tilde{K}(x)\vi(x)\, \vec{\nu}(x)\, d\sigma(x)\right)\\
& = &
\lim_{\eps\to 0}
\left(\int_{|x|\geq\eps}\nabla \tilde{K}(x)\, \vi(x)\, dx +
\int_{|x|=\eps}\tilde{K}(x)\vi(x)\, \frac{x}{|x|}\, d\sigma(x)\right)\, .
\end{eqnarray*}
It remains to prove that
$$
\lim_{\eps\to 0} \int_{|x|=\eps}\tilde{K}(x)\vi(x)\, \frac{x}{|x|}\, d\sigma(x)=
\vi(0)\int_{|x|=1} \tilde{K}(x)\, \frac{x}{|x|}\, d\sigma(x)\, .
$$
Let
$$
c= \int_{|x|=1}\tilde{K}(x)\, \frac{x}{|x|}\, d\sigma(x) =
\int_{|x|=\eps} \tilde{K}(x)\, \frac{x}{|x|}\, d\sigma(x)\, .
$$
The last equality follows from a simple change of variables and homogeneity of
$\tilde{K}$. We have
\begin{eqnarray}
\label{BHeq9}
\lefteqn{\int_{|x|=\eps} \tilde{K}(x)\vi(x)\, \frac{x}{|x|}\, d\sigma(x)}\\
& = &
c\vi(0) +
\int_{|x|=\eps}\tilde{K}(x)\big(\vi(x)-\vi(0)\big)\, \frac{x}{|x|}\, d\sigma(x)\nonumber\\
& \to &
c\vi(0)\nonumber
\end{eqnarray}
as $\eps\to 0$. Indeed,  for $|x|=\eps$
$$
\left|\tilde{K}(x)\big(\vi(x)-\vi(0)\big)\, \frac{x}{|x|}\right|
\leq
C \eps^{1-n}\eps = C\eps^{2-n}\, .
$$
Since the surface area of the sphere $\{|x|=\eps\}$ is $n\omega_n\eps^{n-1}$,
the integral on the right hand side of (\ref{BHeq9}) converges to $0$
as $\eps\to 0$.
\hfill $\Box$

For $\vi\in\S_n$ we have
$$
\nabla (A_\Omega\vi)=c\vi+({\rm p.v.}\, \nabla \tilde{K})*\vi
$$
and hence $\Vert\nabla A_\Omega\vi\Vert_p\leq C\Vert\vi\Vert_p$, $1<p<\infty$,
because the convolution with ${\rm p.v.}\, \nabla \tilde{K}$ is a singular integral.
This and the fact that $A_\Omega:L^p\to L^{p^*}$, $1<p<n$, is bounded yields
$\Vert A_\Omega\vi\Vert_{\dot{W}^{1,p}}\leq C\Vert\vi\Vert_p$, $1<p<n$.
Now the result follows by a density argument.
\hfill $\Box$

\section{Proof of Theorem~\ref{main}}
\label{trzy}

Let
$$
\Phi_t=\tilde{K}\chi_{\bbbr^n\setminus B(0,t)}.
$$
Then
$$
A_\Omega^*f(x)=\sup_{t> 0}|f*\Phi_t(x)|.
$$
Recall that $\Phi_t=(\Phi_t^1,\ldots,\Phi_t^m)$
is a vector valued function and the convolution with $f=(f_1,\ldots,f_m)$ is understood as
$$
f*\Phi_t=\sum_{i=1}^m f_i*\Phi_t^i.
$$
\begin{lemma}
\label{uno}
The distributional gradient of $\Phi_t^i$, $i=1,2,\ldots,m$ equals
$$
\nabla\Phi_t^i=\nabla\tilde{K}_i\chi_{\bbbr^n\setminus B(0,t)}+\mu_t^i,
$$
where $\mu_t^i$ is a measure defined by
$$
\mu_t^i=\frac{1}{t^{n-1}}\Omega_i(x/|x|)\, \frac{x}{|x|}\, \H^{n-1}\rss S^{n-1}(0,t),
$$
i.e.
$$
\int_{\bbbr^n}\vi\, d\mu_t^i=\frac{1}{t^{n-1}}
\int_{S^{n-1}(0,t)} \vi(x)\Omega_i(x/|x|)\, \frac{x}{|x|}\, d\sigma(x)\, .
$$
\end{lemma}
We can also write it in an abbreviated form
$$
\nabla\Phi_t=\nabla\tilde{K}\chi_{\bbbr^n\setminus B(0,t)} + \mu_t,
\quad
\mu_t=(\mu^1_t,\ldots,\mu^m_t).
$$
{\em Proof.}
The lemma is a straightforward consequence of the integration by parts. Indeed, for
$\vi\in C_0^\infty(\bbbr^n)$ we have
\begin{eqnarray*}
\lefteqn{\nabla\Phi_t^i[\vi]
 := 
-\int_{\bbbr^n} \Phi_t^i(x)\nabla\vi(x)\, dx 
=
-\int_{\bbbr^n\setminus B(0,t)} \tilde{K}_i(x)\nabla\vi(x)\, dx} \\
& = &
\int_{\bbbr^n\setminus B(0,t)} \nabla \tilde{K}_i(x)\vi(x)\, dx +
\int_{S^{n-1}(0,t)} \tilde{K}_i(x)\vi(x)\, \frac{x}{|x|}\, d\sigma(x) \\
& = &
\int_{\bbbr^n\setminus B(0,t)} \nabla \tilde{K}_i(x)\vi(x)\, dx +
\frac{1}{t^{n-1}} \int_{S^{n-1}(0,t)} \Omega_i(x/|x|)\vi(x)\, \frac{x}{|x|}\, d\sigma(x).
\end{eqnarray*}
The proof is complete.
\hfill $\Box$

\begin{lemma}
\label{cobra}
The convolution with $\nabla\Phi_t$ is a bounded operator from
$L^p(\bbbr^n,\bbbr^m)$ to $L^p(\bbbr^n)$ for all $1<p<\infty$,
$$
\Vert f*\nabla\Phi_t\Vert_p=
\left\Vert\sum_{i=1}^m f_i*\nabla\Phi_t^i\right\Vert_p
\leq C\Vert f\Vert_p
$$
with the constant $C>0$ independent of $t$.
\end{lemma}
{\em Proof.}
For $f=(f_1,\ldots,f_m)\in L^p$ we have
\begin{eqnarray*}
|f*\nabla\Phi_t|
& \leq &
\sum_{i=1}^m
\left|\int_{|x-z|\geq t} f_i(z)\nabla \tilde{K}_i(x-z)\, dz\right| + |f*\mu_t|   \\
& \leq &
\sum_{i=1}^m \sup_{\tau>0} \left|\int_{|x-z|\geq \tau} f_i(z)\nabla \tilde{K}_i(x-z)\, dz\right| 
+ C|f|*\sigma_t \\
& := &
\sum_{i=1}^m I_i^* f_i + C|f|*\sigma_t,
\end{eqnarray*}
where $\sigma_t$ is the normalized Lebesgue measure on $S^{n-1}(0,t)$
and $I_i^* f_i$ are maximal singular integrals defined independently of $t$.
Now it suffices to observe that the operators $I_i^*$ are bounded in $L^p$ and
$$
\Vert |f|*\sigma_t\Vert_p \leq \Vert f\Vert_p\Vert\sigma_t\Vert = \Vert f\Vert_p
$$
by Young's inequality.
\hfill $\Box$

\begin{lemma}
\label{tre}
The convolution with $\Phi_t$ is a bounded operator from $L^p(\bbbr^n,\bbbr^m)$ to
$\dot{W}^{1,p}(\bbbr^n)$ for $1<p<n$,
$$
\Vert f*\Phi_t\Vert_{\dot{W}^{1,p}} \leq C\Vert f \Vert_p,
\quad
f\in L^p(\bbbr^n,\bbbr^m),
$$
with the constant $C>0$ that does not depend on $t$. Moreover
$$
\nabla(f*\Phi_t) = f*\nabla\Phi_t.
$$
\end{lemma}
{\em Proof.}
Clearly $\Vert f*\Phi_t\Vert_{p^*}\leq C\Vert f\Vert_p$, because $|f*\Phi_t|$ can be estimated by the Riesz potential
of $|f|$.

If $\vi\in C_0^\infty(\bbbr^n,\bbbr^m)$, then $\nabla(\vi*\Phi_t)=\vi*\nabla\Phi_t$ and hence
Lemma~\ref{cobra} yields
$$
\Vert\nabla(\vi*\Phi_t)\Vert_p\leq C\Vert\vi\Vert_p.
$$
The lemma follows now from a density argument.
\hfill $\Box$

\begin{lemma}
\label{due}
The maximal operator
$$
T^*f=\sup_{t>0}|f*\nabla\Phi_t|
$$
is bounded $T^*:L^p(\bbbr^n,\bbbr^m)\to L^p(\bbbr^n)$
for all $n/(n-1)<p<\infty$.
\end{lemma}
{\em Proof.}
Observe that $|f|*\sigma_t\leq \S f$ and hence the estimates from Lemma~\ref{cobra} 
give
$$
T^*f \leq\sum_{i=1}^m I_i^*f_i + C\S f.
$$
Thus the result follows from the boundedness of the maximal singular integrals and
from the boundedness of the spherical maximal function. This is the only moment
where we use the assumption $p>n/(n-1)$.
\hfill $\Box$

We are ready now to complete the proof of the theorem. Clearly for $f\in L^p(\bbbr^n,\bbbr^m)$
we have $|A_\Omega^*f|\leq C\I_1|f|$, so $A_\Omega^*:L^p(\bbbr^n,\bbbr^m)\to L^{p^*}(\bbbr^n)$,
$1<p<n$ is bounded. It suffices to prove now that for $f\in L^p(\bbbr^n,\bbbr^m)$, 
$p>n/(n-1)$, $A_\Omega^*f$ is absolutely continuous on almost all lines
and pointwise partial derivatives of $A_\Omega^*f$ satisfy
\begin{equation}
\label{ac}
|\nabla A_\Omega^* f|\leq T^*f\ \mbox{a.e.}
\end{equation}
Indeed, Lemma~\ref{due} will give the estimate $\Vert \nabla A_\Omega^*f\Vert_p\leq C\Vert f\Vert_p$
for $p>n/(n-1)$ and the result will follow form the characterization of $\dot{W}^{1,p}$ 
by absolute continuity on lines \cite[Section~4.9]{EG}.

To prove absolute continuity along with \eqref{ac} it suffices to show that for almost all 
lines $\ell$ in $\bbbr^n$ the function $A_\Omega^*f$ satisfies
\begin{equation}
\label{ho}
|A_\Omega^*f(x)-A_\Omega^*f(y)|\leq \int_{\overline{xy}} T^*f
\quad
\mbox{for almost all $x,y\in\ell$.}
\end{equation}
Indeed, since the function $T^*f$ is in $L^p$ on almost all lines $\ell$ (Fubini), estimate
\eqref{ho} implies absolute continuity of $A_\Omega^*f$ on compact intervals in $\ell$  with the estimate for
the directional derivatives
\begin{equation}
\label{directional}
D_\nu A_\Omega^*f\leq T^* f
\quad
\mbox{a.e.}
\end{equation}
In particular partial derivatives of $A_\Omega^*f$ are in $L^p$, so $A_\Omega^*f\in \dot{W}^{1,p}$.
Now taking the supremum in \eqref{directional} over a countable and dense set of directions
$\nu\in S^{n-1}$ we obtain
$$
|\nabla A_\Omega^*f|\leq T^*f.
$$
Thus we are left with the proof of \eqref{ho}.

According to Lemma~\ref{tre}, the functions $f*\Phi_t$ are in $\dot{W}^{1,p}$ and hence they are absolutely continuous
on almost all lines. 
Thus for all rational $t>0$ and almost all lines $\ell$, all of the functions $f*\Phi_t$
are absolutely continuous on $\ell$.
Thus for $x,y\in\ell$ we have 
\begin{equation}
\label{all-rational}
\frac{d}{d\tau}(f*\Phi_t)(x+\tau(y-x))=(f*\nabla\Phi_t)(x+\tau(y-x))\cdot(y-x)
\end{equation}
for a.e. $\tau\in (0,1)$ and all rational $t>0$.

Fix $x,y\in\bbbr^n$ such that
$A_\Omega^*f(x)<\infty$, $A_\Omega^*f(y)<\infty$ and  
\eqref{all-rational} is true for all positive rational numbers $t$.

By symmetry we can assume that 
$A_\Omega^*f(x)\geq A_\Omega^*f(y)$. Choose a sequence of positive rational numbers $t_k>0$ such that
$$
A_\Omega^*f(x)=\lim_{k\to\infty} |f*\Phi_{t_k}(x)|.
$$
Clearly
$$
A_\Omega^*f(y)\geq |f*\Phi_{t_k}(y)|
\quad
\mbox{for all $k$.}
$$
Thus
\begin{eqnarray*}
\lefteqn{|A_\Omega^* f(x)-A_\Omega^* f(y)|
 = 
A_\Omega^* f(x)-A_\Omega^* f(y)}\\
& = &
\left( A_\Omega^*f(x)-|f*\Phi_{t_k}(y)|\right) +
\left( |f*\Phi_{t_k}(y)|-A_\Omega^*f(y)\right) \\
& \leq &
A_\Omega^*f(x)-|f*\Phi_{t_k}(y)| \\
& = &
\left( A_\Omega^*f(x)-|f*\Phi_{t_k}(x)|\right)+
\left(|f*\Phi_{t_k}(x)|-|f*\Phi_{t_k}(y)|\right)\, .
\end{eqnarray*}
Passing to the limit yields
$$
|A_\Omega^* f(x)-A_\Omega^* f(y)| \leq
\limsup_{k\to\infty}
|f*\Phi_{t_k}(x)|-|f*\Phi_{t_k}(y)|.
$$
We have
\begin{eqnarray*}
\lefteqn{|f*\Phi_{t_k}(x)|-|f*\Phi_{t_k}(y)|
 \leq 
|f*\Phi_{t_k}(x)-f*\Phi_{t_k}(y)|} \\
& = &
\left|\int_0^1\frac{d}{d\tau}(f*\Phi_{t_k})(x+\tau(y-x))\, d\tau\right| \\
& \leq &
|y-x|\left| \int_0^1 (f*\nabla\Phi_{t_k})(x+\tau(y-x))\, d\tau\right| \\
& \leq &
|y-x|\int_0^1 T^*f(x+\tau(y-x))\, d\tau = \int_{\overline{xy}} T^*f.
\end{eqnarray*}
The proof is complete.
\hfill $\Box$

\section{Proof of Proposition~\ref{spherical}}
\label{cztery}
If $f\in W^{1,1}_{\rm loc}(\bbbr^n)$, then by the trace theorem \cite[Section~4.3]{EG}, for every
$x\in\bbbr^n$ and $t>0$, the integral
$$
\int_{S(x,t)}f(z)\, d\sigma(z)
$$
is well defined and finite.
\begin{lemma}
\label{representation}
If $f\in W^{1,p}$, $1\leq p<n$, then for every $x\in\bbbr^n$ and $t>0$ we have
$$
\barint_{S(x,t)} f(z)\, d\sigma(z) = 
\frac{1}{n\omega_n} \int_{|x-z|\geq t} \nabla f(z)\, \frac{x-z}{|x-z|^n}\, dz.
$$
\end{lemma}
{\em Proof.}
For $\vi\in C_0^\infty(\bbbr^n)$ we have
\begin{eqnarray*}
\barint_{S^{n-1}(x,t)}\vi(z)\, d\sigma(z) 
& = & 
-\int_t^\infty\left(\frac{d}{d\tau}\barint_{S^{n-1}(x,\tau)}\vi(z)\, d\sigma(z)\right)\, d\tau \\
& = &
-\int_t^\infty\barint_{S^{n-1}(x,\tau)}\nabla\vi(z)\, \frac{z-x}{|z-x|}\, d\sigma(z)\, d\tau \\
& = &
\frac{1}{n\omega_n} \int_t^\infty \int_{S^{n-1}(x,\tau)} \nabla\vi(z)\, \frac{x-z}{|x-z|^n}\, d\sigma(z)\, d\tau\\
& = &
\frac{1}{n\omega_n}\int_{|x-z|\geq t}\nabla\vi(z)\, \frac{x-z}{|x-z|^n}\, dz.
\end{eqnarray*}
The case of general $f\in W^{1,p}$ follows by the approximation argument with the use of the H\"older
inequality: $\nabla f\in L^p$, $1\leq p<n$ and the function $z/|z|^n\chi_{\bbbr^n\setminus B(0,t)}$
belongs to $L^{p'}$.
The proof is complete.
\hfill $\Box$

According to the lemma the spherical maximal function of $f\in W^{1,p}$, $1\leq p<n$
can be represented as
$$
\S f(x) =\frac{1}{n\omega_n}\sup_{t>0}\left|\int_{|x-z|\geq t} \nabla f(z)\, \frac{x-z}{|x-z|^n}\, dz\right| =
\frac{1}{n\omega_n}A_\Omega^*(\nabla f)(x)
$$
where $\Omega(z)=z$ for $z\in S^{n-1}$. Thus if the claim of Theorem~\ref{main} remains true for $1<p\leq n/(n-1)$,
then we immediately obtain that the spherical maximal function $\S:W^{1,p}\to\dot{W}^{1,p}$ is bounded for
$1<p\leq n/(n-1)$.
\hfill $\Box$

\end{document}